\documentclass[12pt]{article}

\usepackage{amssymb}
\usepackage[dvips]{graphicx}
\usepackage[english]{babel}
\usepackage{amsfonts,amssymb}

\def\ltitle#1#2{
\begin{center}\bf
#1
\end{center}
\begin{center}\bf
#2
\end{center}
}

\begin{document}

\ltitle{Essential tori in link complements: detecting the satellite structure by monotonic
simplification} {A. Kazantsev\footnote{This work has been partially supported by the grant
NSh-5413.2010.1 -- ``Leading Scientific Schools'' of Russian Federation President's Council for
grants.}}

\par{\bf Abstract.} In a recent work ``Arc-presentation of links: Monotonic simplification''
Ivan Dynnikov showed that each rectangular diagram of the unknot, composite link, or split link can
be monotonically simplified into a trivial, composite, or split diagram, respectively. The
following natural question arises: Is it always possible to simplify monotonically a rectangular
diagram of a satellite knot or link into one where the satellite structure is seen? Here we give a
negative answer to that question both for knot and link cases.\newline

\begin{center}
\par{\bf Introduction}
\end{center}

In \cite{dyn} Ivan Dynnikov shows that by a sequence of elementary moves that do not increase the
complexity of a diagram (number of horizontal or vertical edges) it is possible to transform an
arbitrary rectangular diagram of the unknot into a diagram on which it is obviously seen that the
knot is trivial (See Fig. 1). He also shows that the composite structure of a knot (See Fig. 2) and
split structure of a link (See Fig. 3) can also be found by using elementary moves that do not
increase complexity (further we shall call it \emph{monotonic simplification of a rectangular
diagram}). So, the following question naturally arises: Is it always possible to simplify the
diagram of a satellite knot monotonically to obtain a diagram on which the satellite structure is
seen? Or even more generally: is it always possible to simplify an arbitrary diagram of a link
monotonically to obtain a diagram on which the structure of every incompressible torus~(See Fig. 4)
from the link complement is seen?

\begin{center}
\hfil\includegraphics{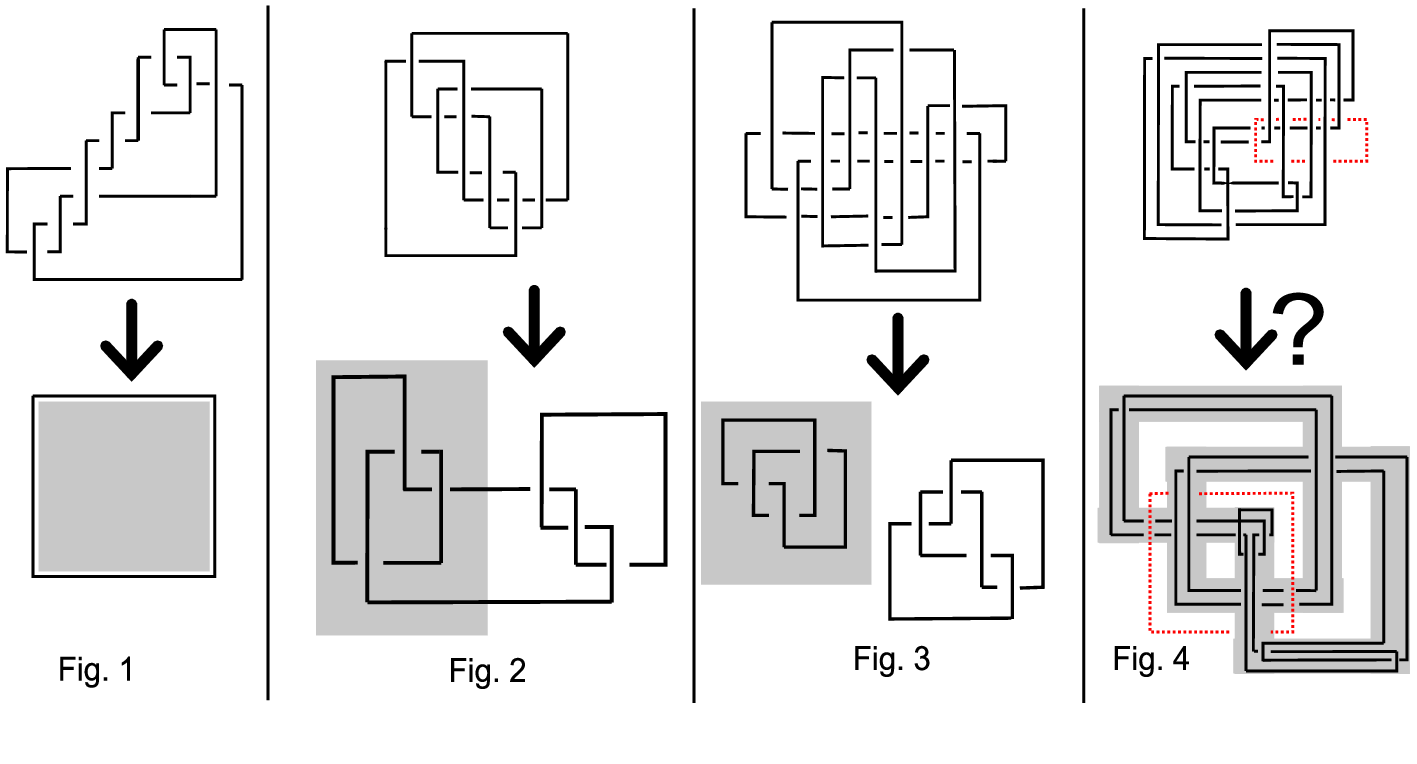} \hfil
\end{center}

The method that we are going to explore originates from works by J.Birman and W.Menasco published
in the beginning of the 90-s in which they developed a technique for studying links and knots
presented in the form of a closed braid (\cite{Bir_IV}, \cite{Bir_V}). The technique uses foliated
surfaces in the complements of links or knots and contains certain tricks which are similar to
those introduced earlier by D. Bennequin in \cite{Ben}.

Using and developing the above mentioned technique J. Birman and W. Menasco published in 1994 a
paper \cite{Bir_SP} where they introduced and studied a class of embedded tori in closed braid
complements. In particular, they proved that each such torus can be admissibly reduced (by the
latter is meant an isotopy of the torus in the complement of a link together with a rearrangement
of the link by exchange moves) so as to have a special position in which the foliation of the torus
induced by open book fibration admits the so-called standard tiling.

However, the geometric description of tori from this class was incomplete (see
\cite{Bir_SP_erratum}). In particular, K. Ng in \cite{Ng} introduced examples of tori admitting a
standard tiling which did not have the form of a thin round tube around a knot, as was stated in
\cite{Bir_SP}. The question of further simplification of such tori into thin tori (we shall use
further this term instead of term ``type k torus'' of~\cite{Bir_SP}) was opened up to this point.
In this paper we will give an answer to an analogous question stated in the language of rectangular
diagrams instead of the closed braids language.

The idea of studying foliated surfaces in link complements was used not only in the case of links
presented in the form of closed braids. In 1995-1996 P. Cromwell (papers \cite{Crom_1995},
\cite{Crom_1996}) adopted this idea for studying so-called arc-presentations of links and
established some of their basic properties. Further development of that technique and results for
rectangular diagrams (which are just another way of thinking about arc-presentations) were obtained
10 years later by I. Dynnikov in \cite{dyn}. In that paper the author shows that a disc, in the
case of a trivial knot, or a sphere in cases of composite and split links can be modified (together
with simplifications of links) so as to have specific foliations induced by open book fibration.
Analogously, one can show that every torus in the complement of an arbitrary arc-presentation can
be modified (together with a simplification of the arc-presentation) to have standard tiling.

So, it is quite natural to restate the question of the further simplification of such tori into
thin tori (whose structure can be seen on the rectangular diagram) in the language of rectangular
diagrams. Exactly this question is investigated in the current article. Examples of
``unsimplifiable'' links and an ``unsimplifiable'' knot rectangular diagrams together with
incompressible non-boundary parallel tori in their complements are presented.

\begin{center}
\par{\bf 1. Preliminaries}
\end{center}

We refer the reader to paper \cite{dyn} for the definitions of \emph{rectangular diagram,
arc-presentation}, and \emph{elementary moves}.

We use the standard terminology in which a knot is a one-component link in $S^3$. Knots and links
are tame and considered up to ambient isotopy.

{\bf Definition.} We shall say that rectangular diagram $K_R$ of a knot $K$ is \emph{a satellite of
rectangular} diagram $L_R$ of link $L$ if the following holds:

1. The x-distance (y-distance) between neighboring vertical (horizontal) edges of $K_R$ is
sufficiently large. By the latter we mean that all distances between neighboring edges are greater
than some constant $\mathrm{Const}$.

2. There is a set $C$ of $L_R$ components whose vertices are lying sufficiently close to vertices
of $K_R$. The vertices of other components are lying sufficiently far from vertices of $K_R$. By
the latter we mean that the distance from every vertex of $C$ and the nearest vertex of $K_R$ is
not greater than, for example $Const/10$. At the same time, there is no component of $L_R$ not
included in $C$ whose some vertex is lying closer than $C/10$ to the vertex of $K_R$.

3. $K_R$ and $L_R$ represent different links.

{\bf Remark.} An example of rectangular diagram and its satellite rectangular diagram can be simply
obtained from Fig. 4 by replacing the ``fat'' shaded knot with its core.

{\bf Definition.} Let $L$ be a link in $S^3$ and $T \hookrightarrow S^3 \backslash L$ is an
essential non-boundary parallel torus. We shall say that the \emph{structure of the torus $T$ is
seen} on a rectangular diagram $L_R$ of link $L$ if $T$ is parallel (in $S^3 \backslash L$) to the
tubular neighborhood of some knot $K$ whose rectangular diagram representation $K_R$ is the
satellite of $L_R$.

{\bf Definition.} Let $L$ be a link in $S^3$ and $T \hookrightarrow S^3 \backslash L$ is an
essential non-boundary parallel torus. We shall say that the \emph{the structure of torus $T$ can
be detected on the rectangular diagram $L_R$ of $L$  by monotonic simplification} if there exist a
sequence of elementary moves not increasing complexity of diagrams which lead from $L_R$ to $L_R'$
where the structure of $T$ is seen.


\begin{center}
\par{\bf 2. Weak example}
\end{center}

{\bf Theorem 1.} There exist rectangular diagrams $L_R$ and $L_R'$ of the link $L$ and an essential
non-boundary parallel torus $T \hookrightarrow S^3 \backslash L$ such that following holds:
\begin{enumerate}
  \item The structure of $T$ is seen on $L_R'$.
  \item The structure of $T$ is not seen on $L_R$.
  \item There is no sequence of elementary moves not increasing
  complexity (destabilizations or exchange moves) converting $L_R$
  into $L_R'$.
\end{enumerate}

\begin{center}
\hfil\includegraphics{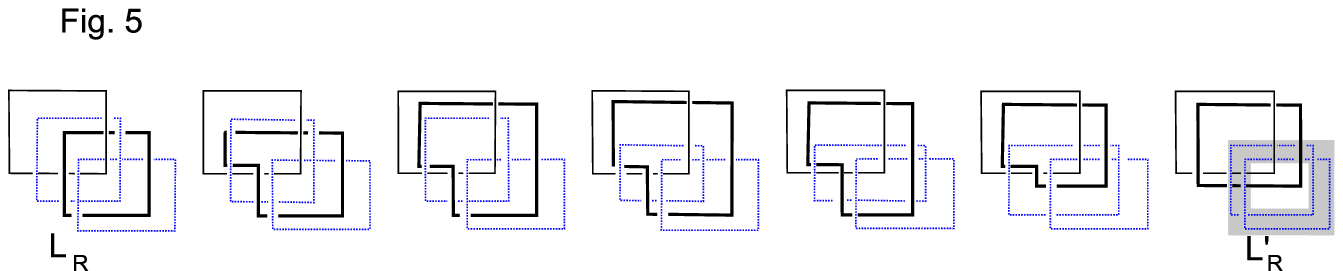} \hfil
\end{center}

{\bf Proof of theorem 1.} An example with diagrams $L_R$ and $L_R'$ is presented on Fig. 5. The
sequence of elementary moves (with increase of complexity!) is shown to verify the fact that two
diagrams represent the same link. It is also obvious that the structure of torus $T$ which includes
two dotted components is seen on $L_R'$ and is not seen on $L_R$. Apart from that, one can easily
verify that there is no possibility to make and exchange move or destabilization with diagram
$L_R$. $\Box$

{\bf Remark.} Actually, the statement of the theorem is not very surprising because of the
following argument. It follows from definition of rectangular diagrams that there can be only
finitely many rectangular diagrams of fixed complexity representing some link. But there exist
links (for example, link $L$ from theorem above) whose complements include infinitely many
non-parallel to each other and boundary non-parallel incompressible tori. Because of that, it is no
wonder that it is impossible to see infinitely many tori on finite number of diagrams.

So, the idea that \textbf{every} incompressible torus from a link complement can be detected by
monotonic simplification on every rectangular diagram is wrong. But may be it is possible to detect
every incompressible torus in the link complement which belongs to some finite set?

This reasoning suggests the following question:

\textbf{Question.} Let $L$ be a link in $S^3$. Then is it true, that every incompressible
non-boundary parallel torus from JSJ-decomposition of $S^3 \backslash L$ can be detected by
monotonic simplification on every rectangular diagram $L_R$ of $L$?

As we shall see further the answer on this question is negative.


\begin{center}
\par{\bf 3. Strong example}
\end{center}

{\bf Theorem 2.} There exist rectangular diagrams $L_R$ and $L_R'$ of a link $L$ and an essential
non-boundary parallel torus $T \hookrightarrow S^3 \backslash L$ \textbf{from JSJ-decomposition} of
$S^3 \backslash L$ such that following holds:
\begin{enumerate}
  \item The structure of $T$ is seen on $L_R'$.
  \item The structure of $T$ is not seen on $L_R$.
  \item There is no sequence of elementary moves not increasing
  complexity (destabilizations or exchange moves) converting $L_R$
  into $L_R'$.
\end{enumerate}

\begin{center}
\hfil\includegraphics{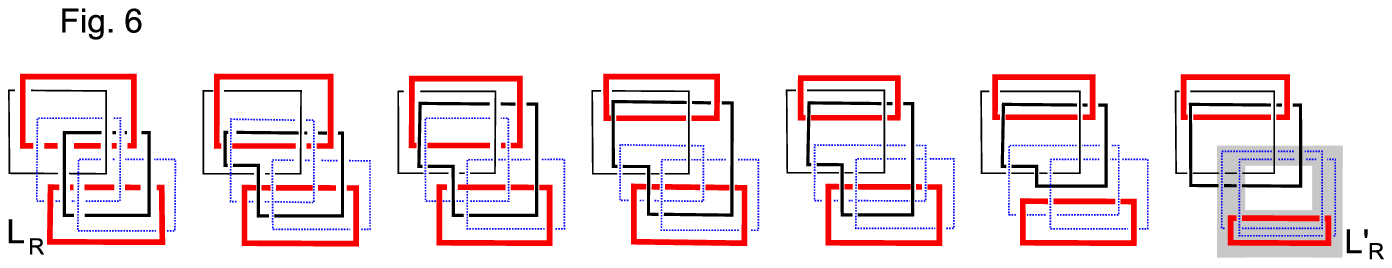} \hfil
\end{center}

{\bf Proof of theorem 2.} An example with diagrams $L_R$ and $L_R'$ is presented on Figure 6. One
can notice that $L_R$ and $L_R'$ differ from corresponding diagrams of theorem 1 only because of
red ''fat'' rings presence. So, the arguments of theorem 1 are absolutely applicable to the current
case.

The only difference from theorem 1 case here is that now the torus (represented by a shaded ring on
Fig. 6) is a torus from JSJ-decomposition.$\Box$

{\bf Remark.} Theorem 2  shows us that a satellite structure cannot always be found by monotonic
simplification when one uses elementary moves presented in the paper~\cite{dyn}. But what if one
will use more a general class of elementary moves? For example, the so-called flypes of rectangular
diagrams, introduced by I. Dynnikov in the paper \cite{Dyn_upsehi}? After all stabilizations,
destabilizations, exchange moves, generalized exchange moves and others are just special cases of
flypes. For this purpose let us recall what a flype is.

{\bf Definition.} Let $R$ be a rectangular diagram and $a$ and $b$ -- positive integers such that
all vertices of $R$ are lying in the following set (it is assumed, that we have coordinate axes):

$\mathbb{R}^2 \backslash (0,a+b) \times (0,a+b)) \cup (0,a) \times (0,b) \cup \{(a+t,t) | t \in
(0,b) \} \cup \{(t,b+t)|t\in(0,a)\}$.

Let us also assume that there is no point from segments $\{a+t,t)|t\in(0,b)\}$ and $\{t,b+t) | t
\in (0,a)\}$, which is not a vertex of a diagram $R$, lying on the intersection of two straight
lines containing edges of diagram $R$.

Then one can construct a new rectangular diagram $R'$, by replacing each vertex $(x,y) \in (0,a)
\times (0,b)$ by vertex $ (y+a,x+b)$. Besides, vertices on the segments $\{(a+t,t) | t\in (0,b)\}$
and $\{(t,b+t) | t \in (0,a)\}$ are deleting or adding in order to obtain the set of vertices of
some rectangular diagram. This can be done uniquely (See Fig. 7 for example).

\begin{center}
\hfil\includegraphics{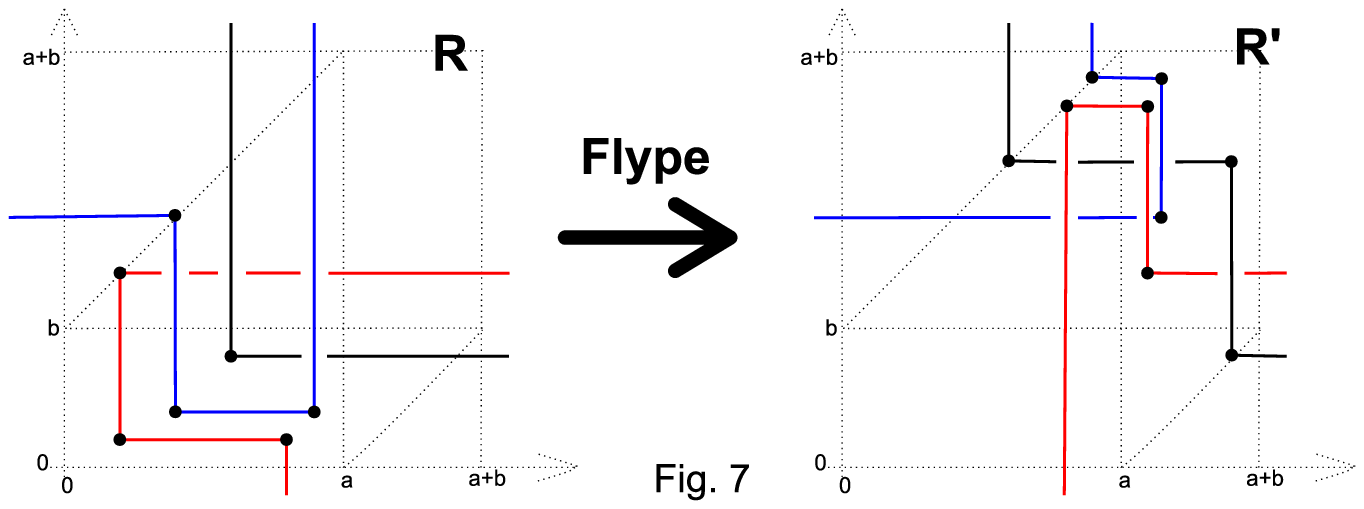} \hfil
\end{center}

The transformation $R \mapsto R'$ is called \emph{flype}. Apart from that every transformation
which is conjugated to the described transformation via a turn on $\pi/2$ or $\pi$ is \emph{also
called flype}.

{\bf Claim.} Links $L_R$ and $L_R'$ from Fig.~6 cannot be obtained from each other by complexity
preserving flypes.

\begin{center}
\par{\bf 4. The relationship of our examples and the example of K. Ng}
\end{center}

Theorem 2 provides us with an example of link $L$ and an incompressible torus $T \in S^3 \backslash
L$ such that $T$ can not be detected by monotonic simplification on rectangular diagram $L_R$ of
$L$.

In the language of arc-presentations, this statement means that there exists an arc-presentation
$L_A$ of link $L$ and an embedding of the torus $T$ into $S^3 \backslash L_A$ such that $L_A$
cannot be modified without increasing complexity, and the torus $T$ cannot be isotoped in the
complement of $L_A$ to be a thin torus (tubular neighborhood of some knot in arc-presentation
form). Let us now describe how the embedding of $T$ can be arranged geometrically.

It turns out that $T$ can be realized in $S^3 \backslash L_A$ so as to have the form of the torus
from the example constructed by K. Ng (Figures 4a, 4b, 4c of paper \cite{Ng}). For clarification of
that fact we provide here a description of an embedding of $L_A$ and $T \in S^3 \backslash \L$ into
$S^3$ by showing an $H_\theta$-sequence of their embedding (Fig. 8) and a foliation induced on $T$
by this embedding (Fig. 9). In our pictures we follow the presentation methods used by K. Ng
(\cite{Ng}, Figures 4, 7).

\begin{center}
\hfil\includegraphics{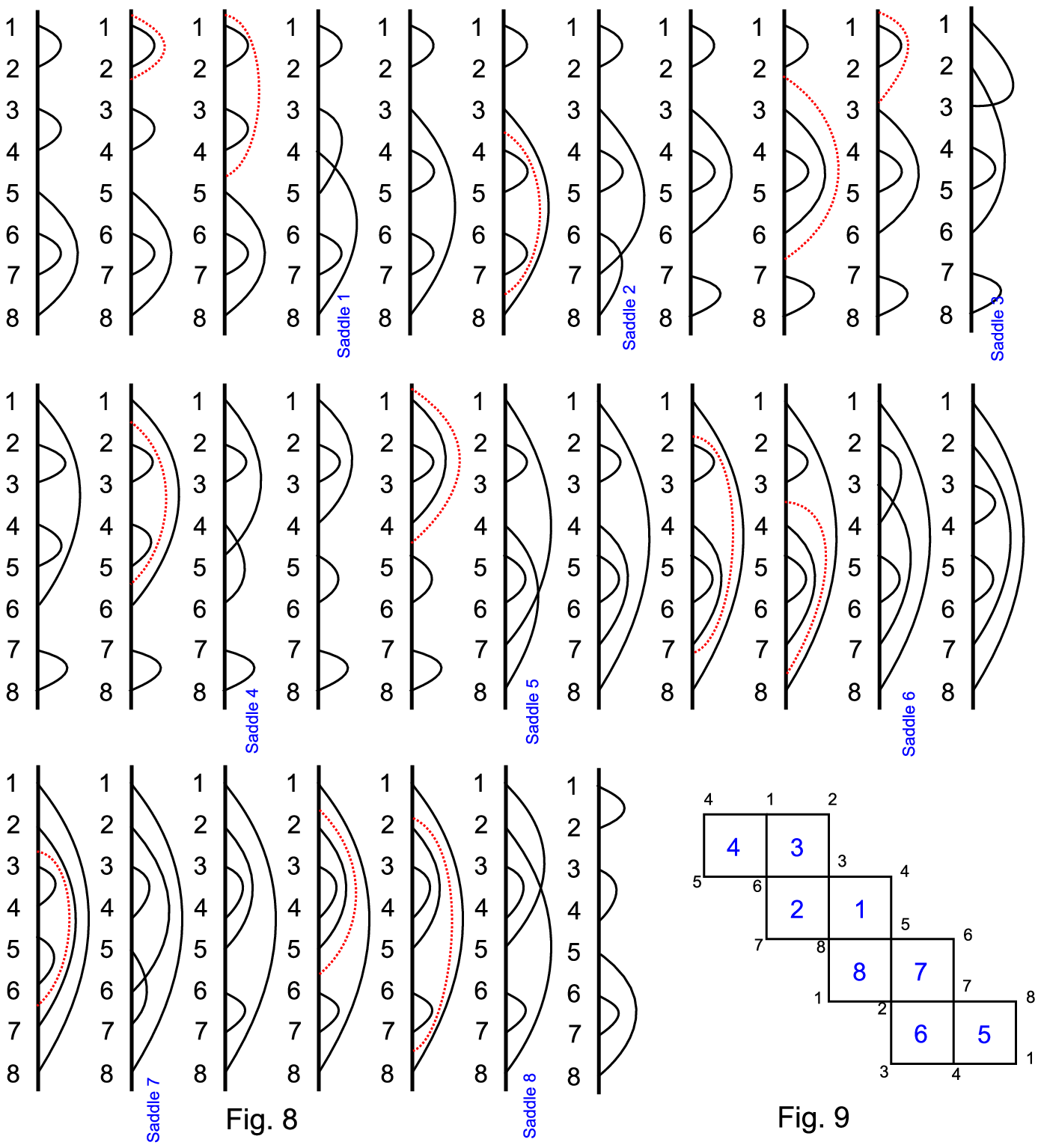} \hfil
\end{center}

Each picture presented in Fig. 8 shows the intersection of $T \cup L$ with corresponding page
$H_\theta$. The first picture corresponds to the case $\theta = 0$, the last picture corresponds to
the case $\theta = 2\pi$ (and that is why these pictures are similar). Black arcs in Fig. 8
represent intersections of torus $T$ with pages $H_\theta$, red dotted arcs represent arcs of the
link $L$. 8 pictures representing critical pages are marked with a legend ``saddle $n$''. Fig. 9
represents standard tiling induced on $T$ by fibration.

\begin{center}
\par{\bf 5. Example with two components}
\end{center}

It was shown that detecting the satellite structure of a link by monotonic simplification of its
rectangular diagram is not possible in general. However, the links in the examples above had four
and six components lying on both sides of the torus, which itself was unknotted. So it is natural
to ask the following question: Can the satellite structure of a satellite link always be detected
by monotonic simplification if the link has fewer than 4 components? We shall answer to this
question in the following theorem:

{\bf Theorem 3.} There exist rectangular diagrams $L_R$ and $L_R'$ of a link $L$ {\bf consisting of
two components} and an essential non-boundary parallel torus $T\hookrightarrow S^3 \backslash L$
such that the following holds:
\begin{enumerate}
  \item The structure of $T$ is seen on $L_R'$.
  \item The structure of $T$ is not seen on $L_R$.
  \item There is no sequence of elementary moves not increasing
  complexity (destabilizations or exchange moves) converting $L_R$
  into $L_R'$.
\end{enumerate}

\begin{center}
\hfil\includegraphics{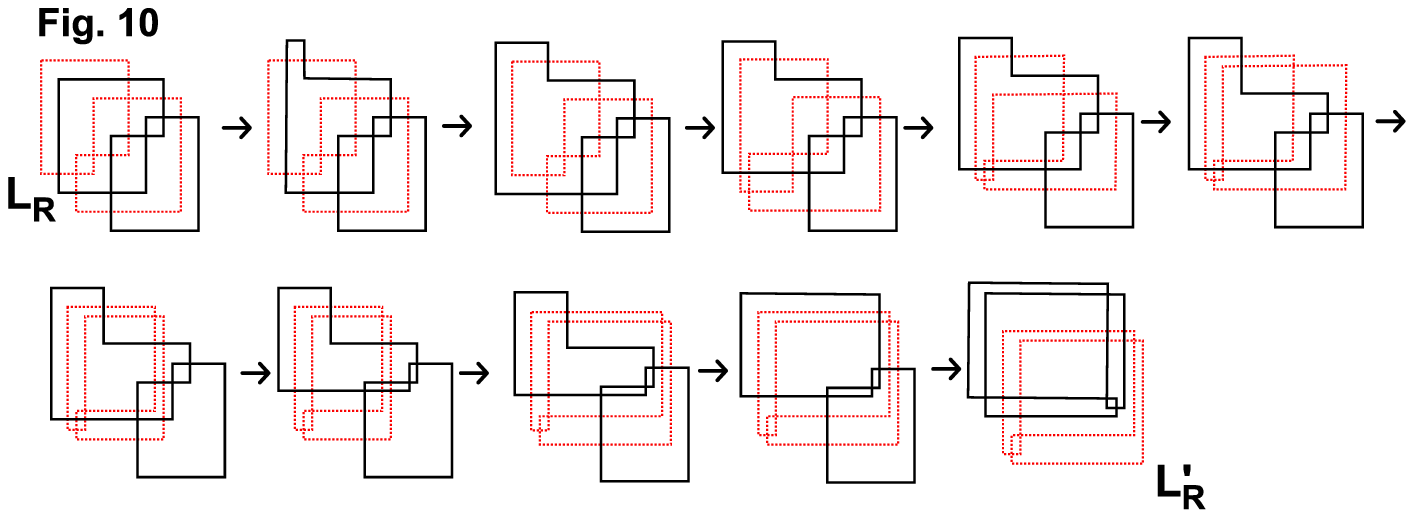} \hfil
\end{center}

{\bf Proof of theorem 3.} An example with diagrams $L_R$ and $L_R'$ is presented on Fig. 10.$\Box$

\begin{center}
\par{\bf 6. Applying foliated surfaces approach to tori.}
\end{center}

It was proved that in case of the link having two or more components lying on both sides of the
torus from the link's complement it is possible to ``lock'' the torus by the link so as the torus
and the link could not be simplified without increase of the complexity. Lets us now pose the
following question:

$\bf{Question.}$ Let $K$ be a {\bf knot} in $S^3$. Is it true, that every incompressible
non-boundary parallel torus from $S^3 \backslash K$ can be detected by monotonic simplification on
every rectangular diagram $L_K$ of $K$?

Using described above standard technique one can show that the torus could be admissibly reduced so
as to have standard tiling. The further simplification is possible if the following holds: The
foliation of the torus has two neighboring (on the binding circle) vertices $a$ and $b$ such that
the arc $ab$ appears when $\theta = t$, exists when $t<\theta<T$ and disappears when $\theta = T$.

$\bf{Definition.}$ Let us call this pattern an \emph{ab-cap} (See Fig. 11). If the segment $ab$ of
binding circle contains no vertices of $K$ (it takes place when, for example, the torus is knotted
and $ab$ lies outside of the torus) then we shall call this ab-cap \emph{empty}.

\begin{center}
\hfil\includegraphics{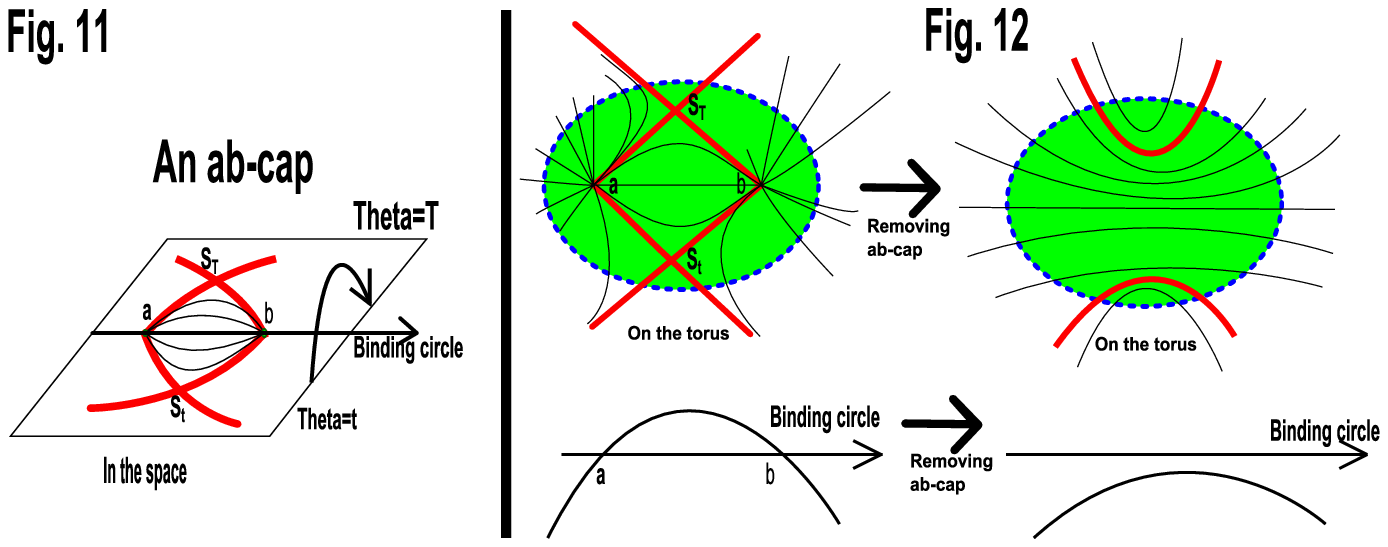} \hfil
\end{center}
$\newline$

If an embedding of a torus has an empty ab-cap then the foliation of the torus induced by open book
fibration will be looking locally as it is presented on Fig. 12 on the top left picture. One can
now see that there is no obstruction to remove the intersection of the ab-cap with the binding
circle (see Fig. 12). After this removal is made the complexity of the torus embedding is decreased
and one can use the standard simplification algorithm leading to new standardly tiled embedding.

{\bf Theorem 4.} If the torus $T$ embedded into the knot complement has the complexity less than 22
then it could be admissibly reduced to thin torus.

{\bf Proof of theorem 4.} Without loss of generality one can suppose that the torus is standardly
tiled. It turns out (this is the result of the computer program written by the author) that if the
standardly tiled torus embedded into $S^3 \setminus K$ has the complexity less than $22$, then
ab-caps exist from both sides of the torus. The torus is embedded into $S^3 \setminus K$, so it has
no knot points inside or outside of it. So, there exists an empty ab-cap, which can be removed.
Therefore the complexity of the embedding will be decreased and one can perform standard
simplification algorithm and receive eventually standardly tiled embedding of smaller complexity.
Performing this procedure consecutively one will eventually receive a thin torus. $\Box$

As we shall see in the next section, there exists an example of the torus of complexity $22$ which
doesn't posses an empty ab-cap.

\begin{center}
\par{\bf 7. An example of a torus and a knot
which could not be simplified}
\end{center}

The conventional way to describe how a torus is embedded into the sphere considers consecutive
intersections of regular and singular pages with the torus. The resulting picture is called an
$H_\theta$-sequence (See. Fig. 8 for example). The same embedding of the torus as in Figure 8 can
be represented by the following picture:

\begin{center}
\hfil\includegraphics{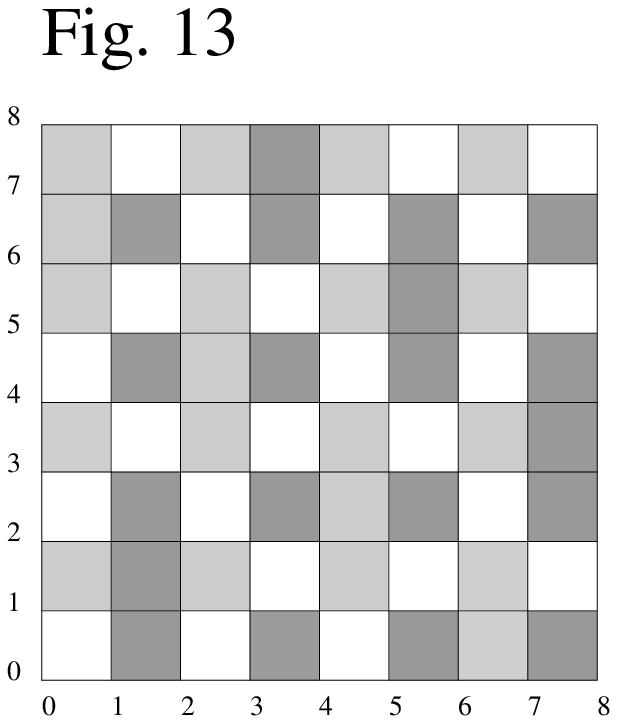} \hfil
\end{center}

This example gives us a new way (introduced by I. Dynnikov) of tori embedding representation. Each
rectangle on the picture above represents the ``life'' of some arc. X-coordinates of horizontal
edges of each rectangle correspond to coordinates of arc endpoints (on the binding circle).
Y-coordinates of vertical rectangle edges represent values of $\theta$ when this arc appears and
disappears. Two rectangles can overlap but no rectangle contains any vertex of another rectangle
(this property comes from the fact that two interleaving arcs cannot exist simultaneously). The
figure above contains two sets of overlapping rectangles (dark grey rectangles of size $1 \times 3$
and light gray of size $3 \times 1$) corresponding to the arcs with odd left point coordinate and
even left point coordinate respectively. As usual, the picture is cyclical. Using this way of
representation we can prove the following theorem.

{\bf Theorem 5.} There exist a knot $K$ and a standardly tiled torus $T \hookrightarrow S^3
\setminus K$ such that there is no empty ab-caps on one side of it.

{\bf Proof of theorem 5.} An example of such torus is presented on the figure below:

\begin{center}
\hfil\includegraphics{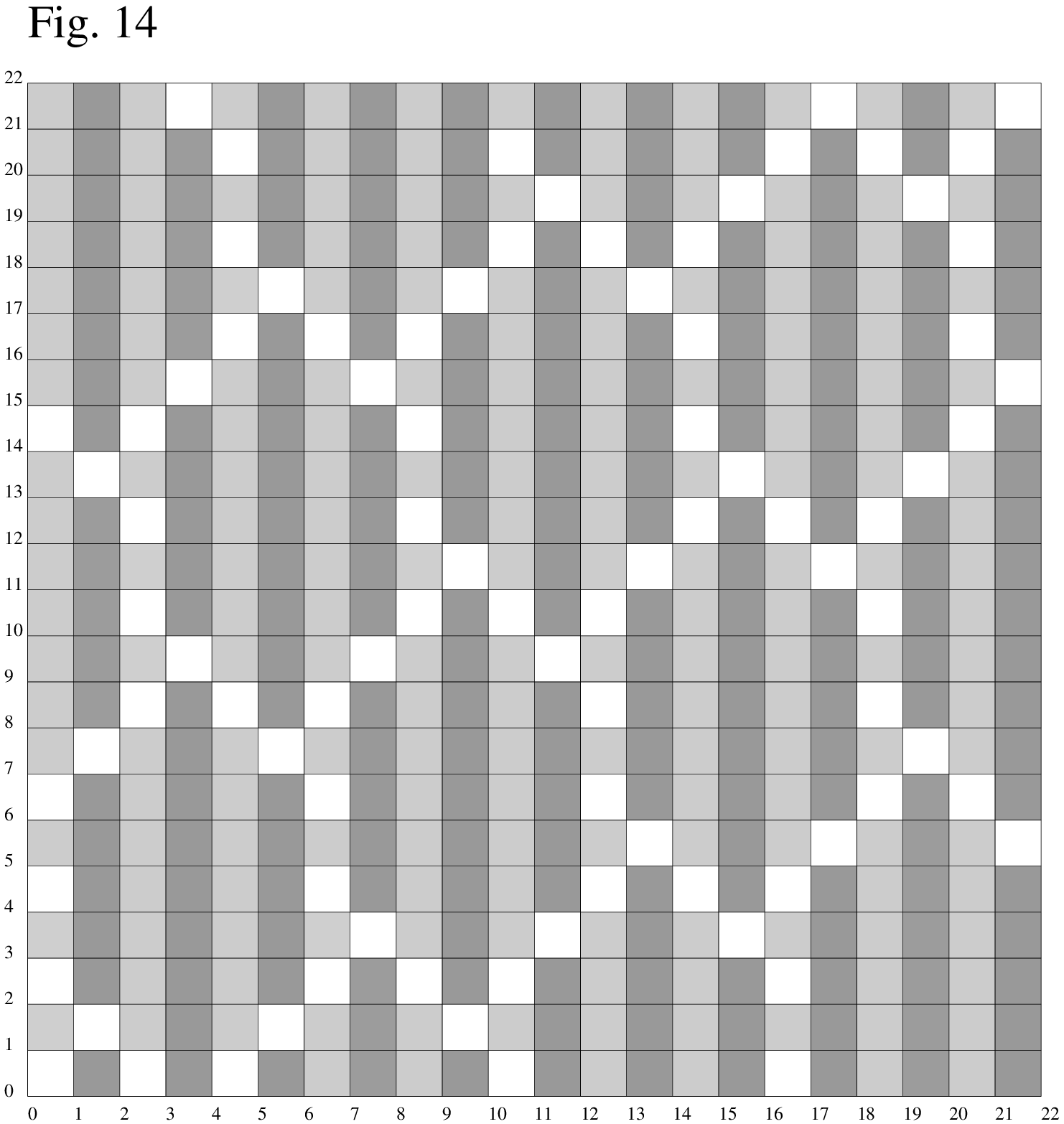} \hfil
\end{center}

One can verify that there are no ab-caps where a is even and b is odd. One also can notice that the
figure determines standardly tiled embedding, because each vertex participates in precisely $4$
saddles. It follows from the fact that there are precisely $4$ rectangles adjacent to each vertical
line -- two from one side and two from another. Apart from that, each saddle rearranges precisely
two arcs into other two arcs (two rectangles are adjacent under and two rectangles are adjacent
over each horizontal line).

In order to explain how the torus is embedded into $S^3$ let us look
at the following figures.

\begin{center}
\hfil\includegraphics{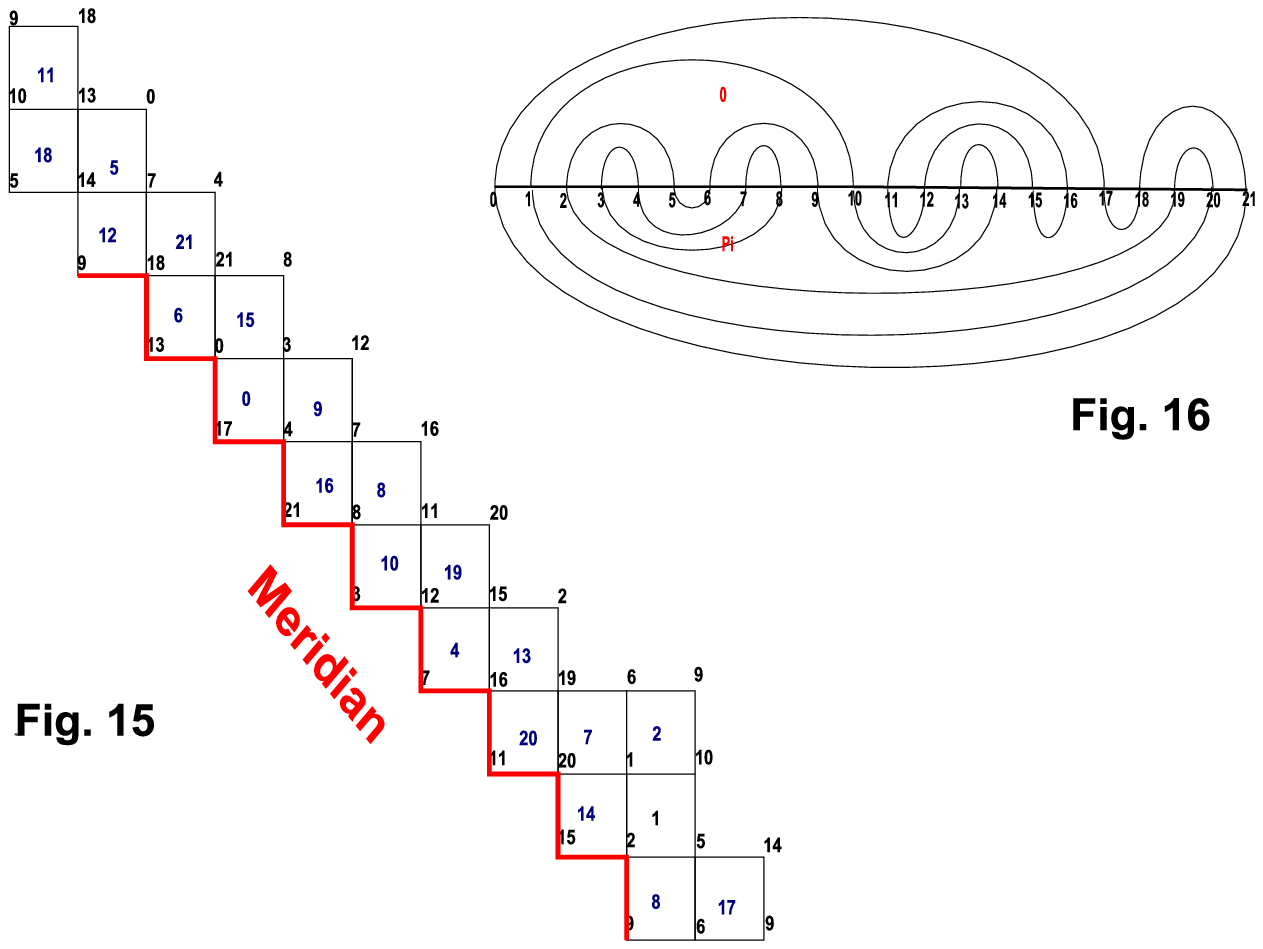} \hfil
\end{center}

Fig. 15 presents a staircase tiling pattern describing the torus in a canonical way. On Fig. 16 the
intersection of the torus with two half-planes corresponding to $\theta = 0$ and $\theta = \pi$ is
pictured. $\Box$

So, we have found a torus $T$, which has no empty ab-caps from one side. It turns out that there is
a way to place a rigid knot into the complement of this torus such that the pair $(T, K)$ could not
be simplified anymore. An example of the knot together with the torus is presented on the following
figure:

\begin{center}
\hfil\includegraphics{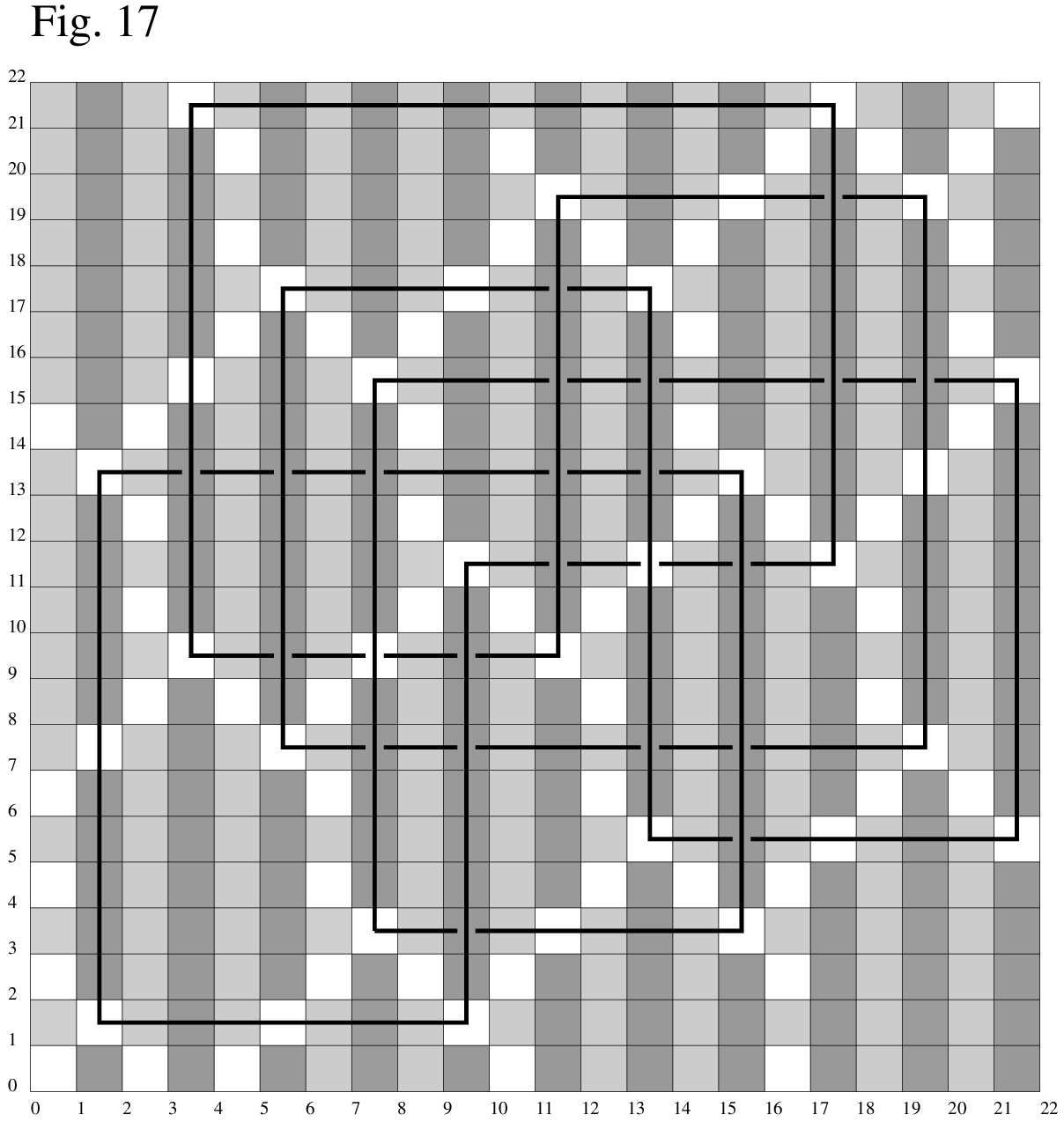} \hfil
\end{center}

The following theorem shows that the knot presented above really gives an answer to the question
stated in paragraph 6.

{\bf Theorem 6.} There exist rectangular diagrams $L_R$ and $L_R'$ of the knot $K$ and an essential
non-boundary parallel torus $T \hookrightarrow S^3 \backslash L$ such that following holds:
\begin{enumerate}
  \item The structure of $T$ is seen on $L_R'$.
  \item The structure of $T$ is not seen on $L_R$.
  \item There is no sequence of elementary moves not increasing
  complexity (destabilizations or exchange moves) converting $L_R$
  into $L_R'$.
\end{enumerate}

{\bf Proof of theorem 6.} An example of diagrams $L_R$ and $L_R'$ is presented on Fig. 18.

\begin{center}
\hfil\includegraphics{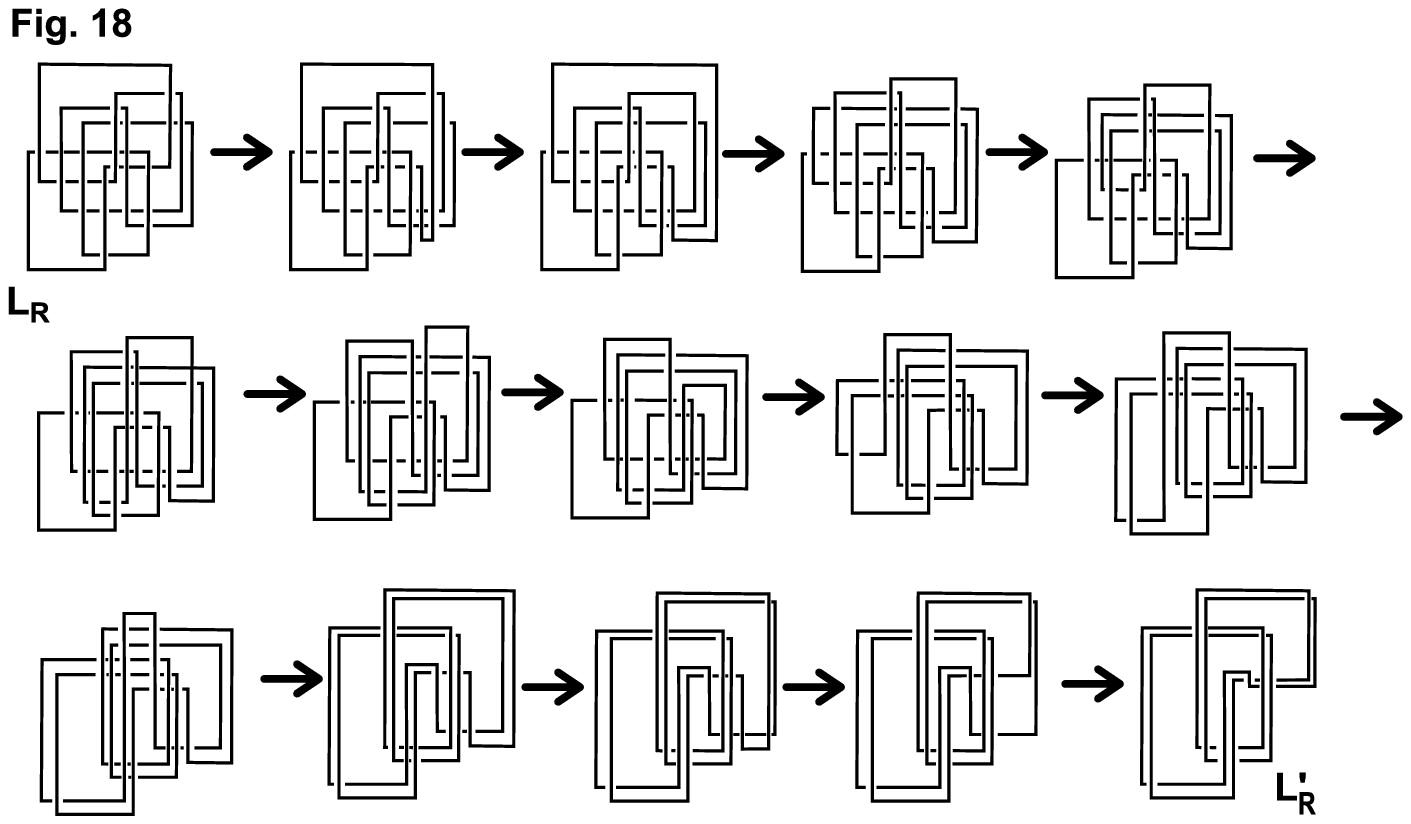} \hfil
\end{center}

$\newline$
{\bf Remark 1.} There are four canonical ways to receive a braid from the rectangular diagram
presented on Fig. 18 (depending on the orientation of the knot and the $90^\circ$-turn).
It is not clear whether one of them gives a counterexample to the initial statement of J. Birman and
W. Menasco about braids (\cite{Bir_SP}, theorem 1) or not.

$\newline$
{\bf Remark 2.} It is possible to make $531$ different flypes without increase of the complexity
with rectangular
diagram presented on Fig. 18 (it is the result of the computer program written by the author).
However,
none of them changes the combinatorial type of the diagram.


\end{document}